\documentclass[11pt]{article}
\usepackage{amsfonts,euscript}
\usepackage{amsmath,amsthm,graphicx}
\usepackage{amssymb,bm,hyperref}
\usepackage{color,pgf,tikz}
\usepackage[a4paper,headsep=0.6cm,left=2cm,right=2cm,top=2cm,bottom=2cm,twoside]{geometry}

\newcommand{\rr}{\mathbb{R}}
\newcommand \nn {\mathbb{N}}

\newcommand \bbc {\mathbb{C}}

 \newcommand {\al} {\alpha}
\newcommand {\be} {\beta}

\newcommand {\Da} {\Delta}

\newcommand {\fy} {\varphi}

\newcommand{\ra}{\rightarrow}
\newcommand{\IN}{{\subset}}

\newcommand {\mmm}{{\setminus}}

\newcommand{\8}{{\infty}}
\newcommand{\0}{{\varnothing}}
\newcommand{\vse}{$\blacksquare$}

\newcommand{\bj}{{\bf {j}}}
\newcommand{\bi}{{\bf {i}}}
\newcommand{\br}{{\bf {r}}}
\newcommand{\bk}{{\bf {k}}}

\newcommand{\bt}{{\bm {t}}}

\newcommand{\eS}{{\EuScript S}}

\newcommand{\eM}{{\EuScript M}}

\newcommand{\eF}{{\EuScript F}}

\def \diam {\mathop{\rm diam}\nolimits}

\def \Lip {\mathop{\rm Lip}\nolimits}
\def \Id {\mathop{\rm Id}\nolimits}

\newtheorem{thm}{\bf Theorem}
 \newtheorem{cor}[thm]{\bf Corollary}
 \newtheorem{lem}[thm]{\bf Lemma}
 \newtheorem{prop}[thm]{\bf Proposition}
 
\theoremstyle{definition}
\newtheorem{ex}{Example}
 
 \newcommand{\beq}{\begin{equation}}
\newcommand{\eeq}{\end{equation}}
\newcommand{\bex}{\begin{ex}}
\newcommand{\eex}{\end{ex}}

\begin{document}

\title{General position theorem and its applications}
\date{}

\author{Vladislav Aseev \and Kirill Kamalutdinov \and Andrey Tetenov \footnote{Supported by Russian Foundation of Basic Research projects 18-01-00420  and 18-501-51021}}

\maketitle

\section*{Introduction}

Consider the following problem:

Let $K$ be the attractor of a system $\eS=\{S_1,\dots,S_m\}$ of contraction maps in $\rr^n$, and let $\dim_H K < n/2$. Suppose that the intersection $S_i(K)\cap S_j(K)$ is nonempty for some $i,j$. Is it possible to change  the maps $S_k\in\eS$ slightly to  maps $S'_k$ to get a system $\eS'=\{S_1',\dots,S_m'\}$ with the attractor $K'$, such that the set $S_i'(K')\cap S_j'(K')$ is empty?\\
To find the answer to this question, we consider
the system $\eS=\eS_0$ as an element of a parametrized family $\eS_t=\{S_{1,t},\dots,S_{m,t}\}$, where the parameter $t$ assumes the values from some subset $D$ in $\rr^n$. We denote the attractor of the system $\eS_t$ by $K_t$. We search for the conditions under which  $S_{i,t}(K_t)\cap S_{j,t}(K_t)$ is empty for almost all $t\in D$. In this case we say  that $S_{i,t}(K_t)$ and $S_{j,t}(K_t)$ are disjoint in general position. 

Particularly, this occurs when Hausdorff dimension of the set $\Da=\{t\in D: S_{i,t}(K_t)\cap S_{j,t}(K_t)\neq\8\}$  is less than  $\dim_H(D)$.

It is possible to make an estimate of $\dim_H(\Da)$ in terms of upper bound for similarity dimensions of the systems $\{\eS_t:\ t\in D\}$.
The method for finding such estimates is based on   General Position Theorem \cite{KT}, which was initially introduced in \cite{TKV}.\\

\section{Definitions and notations}

Let $(X, d)$ be a complete metric space. A mapping $S:\ X \to X$ is a contraction if $\Lip S < 1$ and it  is called a similarity if $d(S(x), S(y)) = r d(x, y)$ for all $x, y\in X$ and some fixed $r$.

Let $\eS=\{S_1,\dots, S_m\}$ be a system of  contractions in a complete metric space $(X, d)$. A nonempty compact set $K\subset X$ is called the attractor of the system $\eS$, if $K = \bigcup \limits_{i = 1}^m S_i (K)$. By Hutchinson's Theorem \cite{Hut}, the attractor $K$ is uniquely defined by the system $\eS$. We also call the set $K$ {\em self-similar} with respect to $\eS$, when all $S_i$ are similarities.

{\bf Multiindices.} Given a system $\eS=\{S_1,\dots,S_m\}$, $I=\{1,\dots,m\}$ is the set of indices, $I^{*}=\bigcup\limits_{n=1}^\infty I^n$ 
is the set of all finite $I$-tuples, or multiindices ${\bf j}=j_1j_2...j_n$. By ${\bf i}{\bf j}$ we denote the concatenation of the corresponding multiindices; we write ${\bf i}\sqsubset{\bf j}$, if  ${\bf j}={\bf i}{\bf k}$ for some ${\bf k}\in I^{*}$; we say that ${\bf i}$ and ${\bf j}$ are {\em incomparable}, if neither ${\bf i}\sqsubset{\bf j}$ nor ${\bf j}\sqsubset{\bf j}$;
by $\bi\wedge\bj$ we mean the maximal $\bk$ for which $\bk\sqsubset\bi$ and $\bk\sqsubset\bj$; by |\bi| we denote the length of $\bi$.

We write $S_{\bf j}=S_{j_1 j_2 \dots j_n}=S_{j_1}S_{j_2}\dots S_{j_n}$ and for the set $A\subset X$ we denote $S_{\bf j}(A)$ by $A_{\bf j}$;  given a set of m ratios $\{r_k, k\in I\}$ we write $r_\bj=r_{j_1}r_{j_2}\dots r_{j_n}$.

  {\bf The Index Space.} $I^{\infty}=\{{\bf i}=i_1 i_2\ldots:\ \ i_k\in I\}$ is the index space;  $\pi:I^{\infty}\rightarrow K$ is the {\em index map}, which sends ${\bf i}\in I^{\infty}$ to  the point $\bigcap\limits_{n=1}^\infty K_{i_1\ldots i_n}$. For a given vector $\br=(r_1,...,r_m)\in(0,1)^m$ we define a metrics $\rho_\br$  on $I^\8$ by $\rho_\br({\bm\al},{\bm\be})=r_{\bm\al\wedge\bm\be}$. The set $I^\8$ supplied with this metrics  will be denoted by $I_{\rho_\br}^\8$. Let $s_\br$ denote the unique solution of the Moran  equation $r_1^s+\dots+r_m^s=1$.  Then, by \cite[Theorem 6.4.3]{Edgar}, $\dim_HI_{\rho_\br}^\8=s_\br$. \\

{\bf Separation conditions.} Denote $\eF=\{S_\bi^{-1} S_\bj: \bi,\bj\in I^{*}\}$. Then the system $\eS=\{S_1,\dots, S_m\}$ of contraction similarities has the Weak Separation Property (WSP) iff $\rm Id \notin \overline{\eF \setminus \Id}$ \cite{Zer}. The system $\eS$ satisfies Open Set Condition (OSC) if there is an open set $V$ such that for any $i\in I$, $S_i(V)\in V$ and for any non-equal $i,j\in I$, $S_i\cap S_j(V)=\0$. The system satisfies Strong Separation Condition (SSC), if for any non-equal $i,j\in I$, $K_i\cap K_j=\0$. There are well-known implications (SSC)$\ra$(OSC)  and  (OSC)$\ra$(WSP) \cite{SSS7,Lau,Zer}\\

\section{General position theorem}

We begin with a simple example. Let $A$, $B$ be compact subsets in $\rr^n$, and the set $B$ is being translated by a vector $t\in D$, where $D\IN \rr^n$. We wish to understand, how large can be the set of  parameters $\Da=\{t\in D:\ A\cap (B+t)\neq\0\}$, which we will call the set of exceptional parameters. 

It's easy to see that $A\cap (B+t)\neq \0$ is equivalent to: " there are such $a\in A$, $b\in B$ that $a=b+t$". Finding  $t$ from this equation, we see that $\Da = \{a-b:\ a\in A,\ b\in B\}$. How to evaluate the Hausdorff dimension of the set $\Da$ in terms of $A$ and $B$?

For that reason we introduce the map $f:\ A\times B\to \Da$, $f(a,b)= a-b$. Since $f$ is Lipschitz,   $\dim_H \Da \le \dim_H (A\times B)$, and if the product $A\times B$ has the dimension less than $\dim_H D$, then $A$ and $B+t$ are disjoint for almost all $t\in D$. 

We will extend this approach to a very general situation, taking a normed linear space $\eM$ instead of $\rr^n$, replacing $A$ and $B$ by metric spaces $(L_1, \sigma_1)$, $(L_2, \sigma_2)$  and finding the set $\Da$ for parametrized families $A_t=\fy_1 (t,L_1)$ and $B_t=\fy_2 (t,L_2)$ instead of $A$ and $B+t$. \cite{KT}:

\begin{thm}\label{genpos}
Let the Cartesian products of metric spaces $(D, \rho)$, $(L_1, \sigma_1)$, $(L_2, \sigma_2)$  be supplied with the canonical metrisation (see \cite[\S 21.VI, (1)]{Kur}). 
Let continuous maps $\varphi_1: D\times L_1 \to {\mathcal M}$ and $\varphi_2: D\times L_2 \to {\mathcal M}$ to the normed linear space $({\mathcal M}, \|.\|)$ be such that:\\
(a)  there are $C_0>0$ and $\alpha>0$ such that for any $i=1,2$ and for all $(\xi, x),(\xi, y)$ in $D\times L_i$  the estimate holds
$$\|\varphi_i(\xi, x)-\varphi_i(\xi, y)\| \leq C_0[\sigma_i(x,y)]^{\alpha}$$
(uniform $\alpha$-H\"older  continuity condition);\\
(b) there are such $M_0>0$ and $\beta>0$ that for any $(x_1,x_2)\in L_1\times L_2$ and $\xi, \xi'\in D$ the function 
$$\Phi(\xi, x_1, x_2): = \varphi_1(\xi, x_1)- \varphi_2(\xi, x_2)$$
on the set $D\times L_1\times L_2$ satisfies the condition 
\begin{equation}\label{antilip }\|\Phi(\xi', x_1,x_2) - \Phi(\xi, x_1, x_2)\| \geq M_0 [\rho(\xi', \xi)]^{\beta}\ .\end{equation}
Then Hausdorff dimension of the set $\Delta: = \{ \xi\in D: \varphi_1(\xi, L_1)\cap \varphi_2(\xi, L_2) \neq \varnothing\}$   satisfies
\begin{equation}\label{dimda }\dim_H \Delta \leq \min \{ (\beta/\alpha)\dim_H(L_1\times L_2), \dim_H D \} \ .\end{equation}
Moreover, if the spaces $(L_1, \sigma_1)$, $(L_2, \sigma_2)$ are compact, $\Da$ is closed in $D$.
\end{thm}
{\bf Proof.}\\ 
Put $\tilde{\Delta}: = \{ (\xi, x_1, x_2)\in D\times L_1\times L_2: \varphi_1(\xi, x_1) = \varphi_2(\xi, x_2)\} =\{ (\xi, x_1, x_2)\in D\times L_1\times L_2: \Phi(\xi, x_1,x_2) =0 \} $ and notice that $\Delta = \mbox{pr}_1 \tilde{\Delta}$, where $\mbox{pr}_1: D\times L_1\times L_2 \to D$ is the canonical projection. 

Applying canonical projection $\mbox{pr}_2: D\times(L_1\times L_2)\to L_1\times L_2$ we obtain a set $\Delta_L : = \mbox{pr}_2 (\tilde{\Delta})$, that is, 
$$ \Delta_L = \{ (x_1, x_2)\in L_1\times L_2 |\ \exists \xi\in D: \varphi_1(\xi, x_1)=\varphi_2(\xi, x_2) \} .$$

The maps $\pi_D = \mbox{pr}_1|_{\tilde{\Delta}}: \tilde{\Delta} \to \Delta$ and $\pi_L = \mbox{pr}_2|_{\tilde{\Delta}}: \tilde{\Delta}\to \Delta_L$ are continuous open maps  (by properties of canonical projections). Let us show that  $\pi_L$ is a bijection. Indeed, if for  $(\xi', x_1', x_2')\in \tilde{\Delta}$ and $(\xi'', x_1'', x_2'')\in \tilde{\Delta}$ the equality $\pi_L(\xi', x_1', x_2') = \pi_L(\xi'', x_1'', x_2'')$ holds, then $(x_1', x_2') = (x_1'', x_2'')= (x_1, x_2)$, whereas $\Phi(\xi', x_1, x_2) = 0 = \Phi(\xi'', x_1, x_2)$. 
Then from (\ref{antilip }) it follows that $0 = \|\Phi(\xi', x_1, x_2) - \Phi(\xi'', x_1, x_2)\| \geq M_0 [\rho(\xi', \xi'')]^{\beta}$, that is, $\rho(\xi', \xi'')=0$. This means that $\xi' = \xi''$. 

Since every open bijective continuous map is a homeomorphism (see \cite[\S 13.XIII]{Kur}), the maps $\pi_L$ and $\pi_L^{-1}$ are homeomorphisms.

Now we find H\"older continuity estimate for a map $g= \pi_D\circ \pi_L^{-1}: \Delta_L \to \Delta$. Let $\xi'= g(x_1', x_2')$ and $\xi = g(x_1, x_2)$. Then $\Phi(\xi', x_1', x_2')= 0 = \Phi(\xi, x_1,x_2)$ and, particularly, $\varphi_1(\xi', x_1') = \varphi_2(\xi', x_2')$. The inequality (\ref{antilip }) gives an estimate 
$$M_0[\rho(\xi', \xi)]^{\beta} \leq \|\Phi(\xi', x_1, x_2) - \Phi(\xi, x_1, x_2)\| = \| \Phi(\xi', x_1, x_2) - 0\| = $$
$$ = \|\varphi_1(\xi', x_1)- \varphi_2(\xi', x_2)\| \leq \|\varphi_1(\xi', x_1)-\varphi_1(\xi', x_1')\| + \|\varphi_1(\xi', x_1')- \varphi_2(\xi', x_2) \| = $$
$$ = \|\varphi_1(\xi', x_1)- \varphi_1(\xi', x_1')\| + \|\varphi_2(\xi', x_2')- \varphi_2(\xi', x_2) \|. $$

Applying the condition (a), we get the inequality
$$ M_0[\rho(\xi', \xi)]^{\beta} \leq C_0[\sigma_1(x_1, x_1')]^{\alpha} + C_0[\sigma_2(x_2, x_2')]^{\alpha}
 \leq 2C_0 \left[\sqrt{\sigma_1(x_1, x_1')^2+ \sigma_2(x_2, x_2')^2}\right]^{\alpha}\ .$$

Denoting by $\tilde{\sigma}$ the metrics of Cartesian product of the spaces $(L_1, \sigma_1)$ and $(L_2, \sigma_2)$, we get  H\"older continuity estimate of the map $g$:
$$\rho(g(x_1', x_2'), g(x_1, x_2)) \leq (2C_0/M_0)^{1/\beta}[\tilde{\sigma}((x_1', x_2'), (x_1, x_2))]^{\alpha/\beta}.$$

Applying \cite[Proposition 2.3]{Fal} and the inequality $\dim_H \Delta_L \leq \dim_H(L_1\times L_2)$, we get the desired relation (\ref{dimda }):
$$\dim_H \Delta = \dim_H g(\Delta_L) \leq (\beta/\alpha) \mbox{dim}_H(L_1\times L_2)\ \ \mbox{and}\ \ \dim_H \Delta \leq \dim_H D  . $$

Since the maps $\fy_i$ are continuous, $\Phi$ is continuous too. The set $\tilde{\Delta}$ is closed in $D\times L_1 \times L_2$ as a set of zeros of $\Phi$.
If $L_1\times L_2$ is compact, the map   $\pi_D$ is proper \cite[Corollary 8 \S 10]{Bur}, therefore the set $\Delta=\pi_D \tilde{\Delta}$ is closed in $D$ (by properties of canonical projections).
$\blacksquare$\\

{\bf Remarks.}\smallskip\\
1. We see from the inequality (\ref{dimda }) that if the product $L_1 \times L_2$ has sufficiently small dimension, then the sets $\fy(t,L_1)$ and $\psi(t,L_2)$ do not intersect for almost all $t\in D$. The proof of the inequality (\ref{dimda }) in the Theorem does not use the condition that the functions $\fy_1$ and $\fy_2$ are continuous with respect to the metrization of product spaces, so this condition may be omitted. It is needed only to show that  $\Da$ is  closed in $D$.\smallskip\\ 
2. The condition (b) in the Theorem may be considered as a form of transversality condition \cite{SimSolUrb}, where $D\subset \rr^n$ is an open set, $\beta=1$ and $\fy_i$ $(i=1,2)$ are the address maps to different copies of a self-similar set, depending of a parameter $\xi\in D$.\smallskip\\
3. Notice that the only information required of the parameter space $D$ is its Hausdorff dimension. Moreover, if $\dim_HD=s$ but the measure $H^s(D)$ is zero, we take some $s'$ satisfying $\dim_H\Da<s'<s$ to see that $\Da$ is negligible in $D$ in a sense that $H^{s'}(D)=\8$ and $H^{s'}(\Da)=0$.\smallskip

For more easy  understanding of the main idea ot the Theorem \ref{genpos} we apply it to
much more simplified settings. Nevertheless even  the following simplified form will be useful for many applications:

\begin{cor}\label{genposS}
Let $A,B,D$ be some subsets of $\rr^n$. 
Let the map  $\fy: D\times B \to \rr^n$  be such that:\\
(a)  there is $C_0>0$  such that  for any $ x, y\in B$ and $t\in D$, 
$\|\varphi(t, x)-\varphi(t, y)\| \leq C_0\|x-y\|$\\
(b) there is such $M_0>0$ that for any $x\in B$ and $t, t'\in D$ \begin{equation}\label{antilipS }\|\fy(t', x) - \fy(t, x )\| \geq M_0 \|t'-t\|\ .\end{equation}
Then Hausdorff dimension of the set $\Da: = \{ t\in D: \varphi(t, B)\cap A \neq \0\}$   satisfies
\begin{equation}\label{dimdaS }\dim_H \Delta \leq \min \{ \dim_H(A\times B), \dim_H D \} \ \end{equation}
Moreover, if $A$ and $B$ are compact and the map $\fy$ is continuous, then $\Da$ is closed in $D$.
\vse
\end{cor}

One can consider several specific applications which may be derived from the Corollary \ref{genposS}:\medskip

\bex {\it If $A,B\subset\bbc$ and $0\notin\bar A$ and $\dim_H A\times B<2$ then for Lebesgue almost all $z\in\bbc$:   $B\cap zA=\0$}.\eex Indeed, let $M_0=\inf\{|z|:z\in A\}$ and for some $C_0> 0$, let $D=\{z: |z|<C_0 \}$. Then the conditions (a) and (b) of the Corollary \ref{genposS}  are fulfilled. Therefore, if $\dim_H(A\times B)<2$ then for Lebesgue almost all $z\in D$ the sets $A$ and $B$ are disjoint. Letting $C_0$ tend to infinity we get that the statement is true for Lebesgue almost all $z\in\bbc$.

\bex {\it If $A,B\IN\rr^n$, $M_2>M_1>0$, a map $f:B\times \rr^n\to\rr^n$ is $M_1$-Lipschitz, and 
$ \dim_H(A\times B)<n$, 
then the set $\Da=\{t\in \rr^n: M_2t+f(B,t)\cap A\neq\0\}$ has zero measure in $\rr^n$ }.\eex
In this case the conditions (a),(b) are fulfilled with $C_0=M_1$ and $M_0=M_2-M_1$. Since the set $\Da$ can be represented also as
$\{t\in \rr^n: f(B,t)\cap M_2t+A\neq\0\}$ this means that if $A$ moves faster that the set $B$ is deformed, for almost all $t$ the set $A$ escapes the intersection with the set $f(B,t)$. 

\bex {\it Suppose $A,B\IN\rr^n$,  a map $F: \rr^n\to\rr^n$ is bi-Lipschitz, and $f:B\times\rr^m$ is defined by $f(x,t)=F(x+t)$.
$ \dim_H(A\times B)<n$, 
then the set $\Da=\{t\in \rr^n: f(B,t)\cap A\neq\0\}$ has zero measure in $\rr^n$ }.\eex

In this case we can interpret $f(B,t)$ as a bi-Lipschitz distortion  of a translation of the set $B$ by a vector $t$.

\section{Application of General Position Theorem to self-similar sets}
The General Position Theorem is a tool for treating more complicated cases, than those in which one of the sets undergoes simple rigid motions or similarities or translations in some curvilinear coordinates. It works with the attractors $K_t$ of parametrized systems $\eS_t$ of contraction maps. These attractors need not be even homeomorphic to each other for different values of the parameter $t$.

To analyze transformations of the attractors of such  systems, we define the following settings for  parametrized families:\medskip\\
{\bf (S1)}. Let $\eS_t =\{S_{1,t},\dots, S_{m,t}\}$ be a system of contraction maps in $\rr^n$, depending on the parameter $t\in D\IN \rr^n$ and let $K_t$ be its attractor.\smallskip\\
{\bf (S2)} Suppose there is a compact set $V$ such that for any $k\in I$ and any $t\in D$, $S_{k,t}(V)\IN V$.\smallskip\\
{\bf (S3)} There is a vector  $\br=(r_1,\dots,r_m)$ such that for any $t\in D$ and for any $k\in I$, $\Lip S_{k_t}\le r_k<1$. Let $\bar r=\max\{r_1,\dots,r_m\}$.\smallskip\\
{\bf (S4)}\label{s4}  There is such $C>0$ that
for any $x\in V$, $k\in I$ and for any $t,t'\in D$, $\|S_{k,t'}(x)-S_{k,t}(x)\|\le C\|t'-t\|$

\subsection{Moving  subpieces apart from each other.}

First notice that it follows from the settings {\bf (S1),(S3)}  that all the address maps are Lipschitz with a constant equal to $\diam(K)$:
 
\begin{lem}\label{param}
If the settings {\bf (S1),(S3)} are fulfilled then the map $\pi: I_{\rho_\br}^\infty \to K$ is $\diam(K)-$Lipschitz.
\end{lem}
{\bf Proof:} (cf. \cite[Ex. 4.2.4]{Edgar}).
Suppose $\bm\al\wedge\bm\be=\bj$, so $\bm\al=\bj\bm\al'$ and $\bm\be=\bj\bm\be'$. From $\rho_\br({\bm\al'},{\bm\be'})=1$ we  get  $\|\pi(\bm\al)-\pi(\bm\be)\|=\|S_\bj(\pi(\bm\al'))-S_\bj(\pi(\bm\be'))|\le r_\bj \diam(K)=\diam(K)\rho_\br({\bm\al},{\bm\be})$. \vse \\

To evaluate the distance between the points in $K_t$ and $K_{t'}$ having the same addresses, we use the Displacement Theorem  for parametrized families ({cf.\cite[Theorem 17]{KT}}): 

\begin{thm}\label{collage}
Suppose  the settings {\bf (S1)---(S4)} hold. Then for any $\bm\alpha\in I^\infty$ and any   $t,t'\in D$ we have \begin{equation}\label{deltapi}\|\pi_{t'} (\bm\alpha)-\pi_t (\bm\alpha)\|\le \dfrac{C\|t'-t\|}{1-\bar r}.\end{equation}
\end{thm}
{\bf Proof:} 
Take $\bm\alpha=i_1 i_2\dots$ and denote $\bm\alpha_k=i_k i_{k+1}\dots$. \\
Since $\pi_t (\bm\alpha_{k})=S_{i_k}^t \pi_t(\bm\alpha_{k+1})$, 
$\|\pi_t(\bm\alpha_{k})-\pi_{t'}(\bm\alpha_{k})\| \le \|S_{i_k}^t\pi_t (\bm\alpha_{k+1})-S_{i_k}^t \pi_{t'}(\bm\alpha_{k+1})\|+ \|S_{i_k}^t \pi_{t'}(\bm\alpha_{k+1})-S_{i_k}^{t'}\pi_{t'}(\bm\alpha_{k+1})\|,$ so
$\|\pi_t(\bm\alpha_{k})-\pi_{t'}(\bm\alpha_{k})\| \le r_{i_k} \|\pi_t(\bm\alpha_{k+1})-\pi_{t'}(\bm\alpha_{k+1})\|+C\|t'-t\|$ for any $k\in \nn$.\\
Therefore $\|\pi_t(\bm\alpha)-\pi_{t'}(\bm\alpha)\|\le {\bar r}^{n+1} \|\pi_t(\bm\alpha_{n+1})-\pi_{t'}(\bm\alpha_{n+1})\|+ C \|t'-t\| \sum\limits_{k=0}^{n} {\bar r}^k$, which becomes (\ref{deltapi}) as $k$ tends to $\8$. 
\vse \\

The following Theorem gives the conditions under which the pieces $K_{\bj,t}$ and $K_{\bk,t}$ are disjoint for almost all $t\in D$:

\begin{thm}\label{main}
Suppose  the settings {\bf (S1)---(S4)} hold. Let $\bj,\bk\in I^*$ be incomparable multiindices.\\
Suppose there are such $c_\bj>0,C_\bk>0$  that for any $x\in V$ and for any $t,t'\in D$,  \beq\|S^{t'}_\bk(x)-S^t_\bk(x)\|\le C_\bk\|t'-t\|\mbox{   and   }\|S_{\bj,t'}(x)-S_{\bj,t}(x)\|\ge c_\bj\|t'-t\|\label{eq6}\eeq
If \beq c_\bj-C_\bk-\dfrac{(r_\bj+r_\bk)C}{1-\bar r}>0\label{eq7}\eeq and  $s_\br<\dim_H(D)/2$, then $K_\bj\cap K_\bk=\0$ for almost all $t\in D$. 
\end{thm}
{\bf Proof:}\quad
Let $\varphi(t,x)=S_{\bk,t} (\pi_t (x))$,\quad $\psi(t,x)=S_{\bj,t} (\pi_t (x))$,\quad $\Phi(t,x,y)=\varphi(t,x)-\psi(t,y)$,\\ $\Delta=\{t\in D:\ K_\bj\cap K_\bk\neq\varnothing\}$.
Note that
\begin{align*}
\|\Phi(t',x,y)-\Phi(t,x,y)\|\ge &\|\psi(t',y)-\psi(t,y)\|-\|\varphi(t',x)-\varphi(t,x)\|;\\
\|\varphi(t',x)-\varphi(t,x)\| \le &\|S_{\bk,t'} (\pi_{t'} (x))-S_{\bk,t} (\pi_{t'} (x))\|+\|S_{\bk,t} (\pi_{t'} (x))-S_{\bk,t} (\pi_t (x))\|;\\
\|\psi(t',x)-\psi(t,x)\| \ge &\|S_{\bj,t'} (\pi_{t'} (x))-S_{\bj,t} (\pi_{t'} (x))\|-\|S_{\bj,t} (\pi_{t'} (x))-S_{\bj,t} (\pi_t (x))\|.
\end{align*}
From Theorem \ref{collage} we have upper estimates\\
$$\|S_{\bk,t} (\pi_{t'} (x))-S_{\bk,t} (\pi_t (x))\| \le \dfrac{r_\bk C\|t'-t\|}{1-\bar r}\mbox{\quad and \quad}
\|S_{\bj,t} (\pi_{t'} (x))-S_{\bj,t} (\pi_t (x))\| \le \dfrac{r_\bj C\|t'-t\|}{1-\bar r}$$
Combining them with inequalities \eqref{eq6}, we obtain
\beq\|\Phi(t',x,y)-\Phi(t,x,y)\| \ge \left( c_\bj - C_\bk - \dfrac{C(r_\bk + r_\bj)}{1-\bar r}\right) \|t'-t\|\label{mlast}\eeq
 Applying the Theorem \ref{genpos} with $\al=\be=1$ we get $\dim_H \Delta<2\dim_H (I^\infty_{\rho_\br})=2s_\br$. \\
Since $s_\br < \dim_H(D)/2$ we get $H^{2s_\br}(\Delta)=0$ and at the same time $H^{2s_\br}(D)=\infty$.
$\blacksquare$

\subsubsection{The case when the parameters are translation vectors.}

We consider the case is when the initial system $\eS=\{S_1,...,S_m\}$ consists of the contraction maps $S_k$ in $\rr^n$; and we consider a parametrized system $\eS_\bt=\{S_{1,\bt},...,S_{m,\bt}\}$ where each $S_{k,\bt}$ is defined by the formula $S_{k,\bt}(x)=S_k(x)+t_k$, where $\bt=(t_1,...,t_m)\in (\rr^n)^m$. Translations have no effect upon the contraction ratios, therefore $\Lip S_{k,t}=r_k$ for any $t$.

First we allow only one map, say $S_{m,t}$, to depend on the parameter $t$, leaving all others unchanged.  

\begin{cor}
Let $\eS_t =\{S_1,\dots, S_{m-1},S_{m,t}(x)=S_m(x)+t\}$ be a system of contraction maps in $\rr^n$, depending on the parameter $t\in \rr^n$ and let $K_t$ be its attractor. Let $1\le k<m$. If  $r_k+r_m+\bar r<1$ and $s_{\br}<n/2$, then $K_{k,t}\cap K_{m,t}=\0$ for almost all $t\in\rr^n$.
\end{cor}
{\bf Proof:}  For any open bounded $D\subset \rr^n$ there is such $V\IN\rr^n$ that  the system $\eS^t $ satisfies the settings {\bf (S1)---(S4)}; since  $C=1$  the condition \ref{eq7} of the Theorem \ref{main} becomes equivalent to $r_k+r_m+\bar r<1$. Therefore $K_{k,t}\cap K_{m,t}=\0$ for almost all $t\in D \subset\rr^n$. The result does not depend on the choice of $D\subset \rr^n$, so it holds for the whole $\rr^n$. $\blacksquare$

Now, if we apply a translation  by  some vector $t_k\in \rr^n$ to each map $S_k\in \eS$, we obtain the following:

\begin{cor}
Let $\eS =\{S_1,\dots, S_{m-1},S_m\}$ be a system of contraction maps in $\rr^n$. Let $\bm t=\{t_1,...,t_m\}$, where $t_k\in \rr^n$. Let $S_{k,\bm t}(x)=S_k(x)+t_k$. Let $K_{\bm t}$ be the attractor of  the system $\eS_{\bm t}=\{S_{1,\bt},...,S_{m,\bt}\}$. If for any non-equal $j,k\in I$, $r_j+r_k+\bar r<1$ and $s_{\br}<n/2$, then for almost all $\bt\in\rr^{mn}$, the system $\eS$ satisfies Strong Separation Condition.\end{cor}
{\bf Proof:} Notice that by Theorem \ref{collage} the maps $\pi_{j,\bt}:I^\8\times\rr^{nm}\to \rr^n$ are continuous with respect to $\bt$. Therefore the function $\rho_{jk}(\bt)=\min\{\|\pi_{j,\bt}(\bm\al)-\pi_{k,\bt}(\bm\be)\|,\bm\al,\bm\be\in I^\8\}$ 
is continuous with respect to $\bt$. Therefore the set $\Da_{jk}=\rho^{-1}(\{0\})$ is closed in $\rr^{nm}$. Since all of its $k$-slices $\{(t_1,..,t_{k-1},t,t_{k+1},...,t_m)\in\Da_{jk};t\in\rr^n\}$ have zero Lebesgue n-dimensional measure, the set $\Da_{jk}$ has zero measure in $\rr^{mn}$. Thus, the set $\Da=\bigcup\limits_{j,k\in I}\Da_{jk}$ also has zero measure in $\rr^{mn}$. Therefore, for almost all $\bt\in \rr^{mn}$, the system $\eS_\bt$ satisfies  Strong Separation Condition.  \vse

\subsection{Non-empty overlaps of prescribed type.}

If we we  get rid of all overlaps in a self-similar set, we obtain a system $\eS$, which satisfy Strong Separation Condition and whose attractor $K$ is just a Cantor set. There is a mush more interesting case, when we use our techniques to obtain a system $\eS$ of contraction maps which has the  attractor $K$ such that the intersections of its pieces $K_j$ strictly follow some predefined pattern. The attractors of such systems possess a set of interesting properties and often they do not satisfy WSP. In this subsection we will see\\ a) how to 
find systems $\eS$ for which two maps $S_1$ and $S_2$ commute and for which $S_1(K)\cap S_2(K)$ is exactly equal to $S_{12}(K)$ and\\
b) how to find systems $\eS$ which do not satisfy OSC though all the pieces $S_i(K)$ are disjoint except $S_1(K)\cap S_2(K)$ which is a single point.

\subsubsection{Exact overlaps: an example}

First we consider the systems $\eS$ in which two maps $S_1,S_2$ have a common fixed point and commute. (cf.\cite{Bar2})
Let the system $\eS_{t}$ in $[0,1]$ consist of  3 maps: $S_1 (x)=tx$, $S_2 (x)=bx$, $S_3 (x)=\dfrac{x+8}{9}$ in $\rr$, where $b,t \in (0,1/9)$. It  depends on the parameter $t$, while $b$ is a fixed value.

Since  the maps $S_{1,t}$ and $S_2$  commute, we have the following inclusion:
\beq S_{1,t}S_2(K_t)\subseteq S_{1,t}(K_t)\cap S_2(K_t) \label{ovl}\eeq
We want to  study for which $t\in(0,1/9)$ the inclusion (\ref{ovl}) becomes equality. In this case we say the system $\eS_t$ has exact overlap $ S_{1}(K)\cap S_2(K)=S_{12}(K)$.

Notice that the same way as in (\cite[Proposition 2(v)]{KT}),
\begin{equation}\label{abc1}K_t\setminus\{0\}=\bigcup\limits_{m,n=0}^\infty S_1^m S_2^n (K_{3,t})\end{equation}

Since $t,b<1/9$ and $K_3\subset [8/9,1]$, for any $m\neq n$, $S_i^m (K_3)\cap S_i^n (K_3)=\0$  for $i=1,2$.

Following the argument of \cite[Proposition 3]{KT} we obtain 
\begin{prop}\label{disjoint}For the system $\eS_{t}$  the following statements are equivalent:\\
(i) For any $m,n\in\mathbb N$, $S_1 ^m(K_3)\cap S_2 ^n(K_3)=\0$;\\
(ii) $K=\{0\}\cup\bigsqcup\limits_{m,n=0}^{\infty}S_1^m S_2^n (K_3)$;\\
(iii) For any $m,n\in\mathbb N$, $S_1^m (K)\cap S_2^n (K)=S_1^m S_2^n (K)$.\vse
\end{prop}

\begin{prop}
The system $\eS_t$ has exact overlap $ S_{1}(K)\cap S_2(K)=S_{12}(K)$  for Lebesgue almost all $t\in (0,1/9)$.
\end{prop}

{\bf Proof:} By Proposition \ref{disjoint} it suffices to find the set of those $t$, for which $S_1 ^m(K_3)\cap S_2 ^n(K_3)=\0$ for any $m\neq n$.\\

Take  non-equal $m,n\in\nn$ and let $D_{mn}=\{t\in (0,1/9):S^m_{1,t}([8/9,1])\cap S^n_2 ([8/9,1])\neq\0\}$.\\

 If $t\in D_{mn}$ then  $\dfrac{8b^n}{9}\le t^m \le\min\left\{\dfrac{9b^n}{8},\dfrac{1}{9^m}\right\}$. Put
$\bar t=\left(\min\left\{\dfrac{9b^n}{8},\dfrac{1}{9^m}\right\}\right)^{1/m}.$\\

To apply the Theorem \ref{main} we interpret the case under consideration in terms of its settings:\\
The system $\eS_{t}$   depends on the parameter $t\in D_{mn}$.\\ The set $V=[0,1]$,  the constant $C=1$. Since the vector $\br=
(\bar t,b,1/9)$, we have $s_\br< 1/2$.\\ Further,
 $S_\bj=S_{1,t}^m$, $S_\bk=S_{2}^n$, therefore $r_\bj={\bar t}^m<\dfrac{9b^n}{8}$, $r_\bk=b^n$.\smallskip\\
 By definition, $c_\bj=\inf\limits_{t,t'\in D_{mn}}\dfrac{{t'}^m-t^m}{t'-t}=\inf\limits_{t\in D_{mn}}mt^{m-1}\ge \inf\limits_{t\in D_{mn}}\dfrac{t^{m}}{t}$.\\ Replacing $t^m$ by $\dfrac{8b^n}{9}$ and  $t$ in denominator by $1/9$, we  get $c_\bj>8b^n$.\smallskip\\
Since $C_\bk=0$, we have $c_\bj-C_\bk-\dfrac{r_\bj+r_\bk}{1-\bar r}>\left(8-\dfrac{9/8+1}{8/9}\right)b^n$.\smallskip\\
Therefore by Theorem \ref{main}, the set $\Da_{mn}=\{t\in D:S^m_{1,t}(K_{3,t})\cap S^n_2 (K_{3,t})\neq\0\}$ is a closed subset of $D_{mn}$ and $\dim_H(\Da_{mn})<1$. 

Let $\Da$ be the union of all $\Da_{mn}$, where $m,n\in\nn$ and $m\neq n$.\\ Then $\dim_H(\Da)\le 2s_\br<1$ which implies the statement of the Proposition.\vse

 For almost all $t$ the systems $\eS_t$  possess several remarkable properties:\\

{\bf 1. Violation of WSP.}  Consider the set $D^*$ of  those values of  the parameter $t\in D\mmm\Da$ for which $\dfrac{\log t}{\log b}$ is irrational. The set $D^*$ has full measure in $D$. For each  $t\in D^*$ there are sequences of positive integers $l_k$,$n_k$ such that the sequence $t^{l_k}b^{-n_k}$ converges to 1. Therefore the system $\eS_{t}$ does not satisfy Weak Separation Property. \\

{\bf 2. Measure and dimension.} The Hausdorff dimension $s$ of the attractor $K_{t}$, $t\in D^*$ is equal to the solution of the equation $t^x+b^x-t^xb^x+9^{-x}=1$. Since the Weak Separation Property is violated, the Hausdorff measure $H^{s}(K_{t_0})=0$.\\

{\bf 3. All $K_t$ are isomorphic.}  For any two sets $K_{t_1}$, $K_{t_2}$, $t_i\in D^*$  there is a homeomorphism $\fy:K_{t_1}\to K_{t_2}$, which agrees with the systems $\eS_1$ and $\eS_2$, i.e. for any $k=1,...,4$ and for any $x\in K_t$, 
$\fy(S_{k,t}(x)= S_{k,t'}(\fy(x))$.\\

We refer the reader to \cite{KT} for  detailed proofs of the  properties of such type of self-similar sets.
 
\subsubsection{One-point intersections: an example}

Take $p,q,r$ in $(0,  1/36)$ and put  $h =\dfrac{1}{2},a=\dfrac{1}{3}$.
Consider a system $\eS=\{S_1,S_2,...,S_6\}$ of contractions in $[0,1]$ whose equations are
$$S_1 (x)=p x,\quad S_2 (x)=a +rx,\quad S_3 (x)=h -qx, \quad S_4 (x)=h-r+rx,$$ $$ \quad S_5 (x)=1-a -rx,\quad S_6 (x)=1-r+rx $$

The similarity dimension for any such system is strictly less than  $1/2$.

Let  $K$ be the attractor of the system $\eS$ and $K_i=S_i(K)$ be its pieces. By the construction, $\{0,1\}\IN K\IN[0,1]$ and the pieces $K_i, i\in \{1,2,3,5,6\}$ are contained in disjoint segments of length $1/36$, while $K_3\cup K_4\IN[h-1/36,h]$ and $K_3\cap K_4\ni\{h\}$ which is the only possible non-empty intersection of the pieces.

We wish to know the set of those $p,q,r$ for which $K_3\cap K_4=\{h\}$. In this case we say that the system $\eS$ has {\em unique one-point intersection}. \\

If $\dfrac{\log p}{\log r}\notin\mathbb Q$, then the system $\eS$  does not have WSP for any $q$. Indeed, consider the maps
$H_m(x)=S_3 S_1^m S_5(x)$ and $G_n(x)=S_4 S_6^n S_2(x)$. Notice that
for any $q>0$ there is a sequence  $(m_k,n_k)\in \mathbb N^2$, such that $ p^{-m_k}r^{n_k+1}$ converges to  $q$ as $k\to\infty$.
Easy computation shows that if we choose such a sequence $(m_k,n_k)$, then the sequence $$G_{n_k}^{-1} H_{m_k} (x) = \dfrac{(r^{n_k+1}-p^{m_k}q)(1-a)}{r^{n_k+2}} + \dfrac{p^{m_k}q}{r^{n_k+1}} x$$ converges to identity,  which means violation of WSP.\\

Therefore we fix  some $p,r\in (0,1/36)$ such that $\log_rp$ is irrational and consider a 1-parameter family of systems $\eS_q$, $q\in (0,1/36)$, for which we show that
 for Lebesgue almost all  $q\in (0,1/36)$ the system $\eS_q$ has unique one-point intersection and does not have
Weak Separation Property.   

For the simplicity of notation, we   denote the system  under consideration by  $\eS$, keeping in mind that it depends on the parameter $q$ whenever  it does not cause  any ambiguity.

From the representation of the pieces $K_3$ and $K_4$ as  unions of infinite sequences 
$$K_3=\{h\}\cup\bigcup\limits_{m=0}^\infty S_3 S_1^m (K \setminus K_1),\ \ K_4=\{h\}\cup\bigcup\limits_{n=0}^\infty S_4 S_6^n (K \setminus K_6),$$
we see that $K_3 \cap K_4=\{h\}$  iff\\ \beq\mbox{ for any }m,n\in \mathbb N\cup\{0\}\mbox{ and any }i\in I\setminus\{6\},\ \  j\in I\setminus \{1\},\quad S_3 S_1^m (K_j)\cap S_4 S_6^n (K_i)= \0 \eeq\\
Note that if $p^m [aq,q] \cap r^{n+1} [a,1] = \0$ then for any $i\in I\mmm\{6\}$, $j\in I\mmm\{1\}$  the intersections $S_3 S_1^m S_j (K) \cap S_4 S_6^n S_i (K)$ are empty. Therefore, in search of those $q$ for which 
$S_3 S_1^m S_j (K)$ and  $S_4 S_6^n S_i (K)$ may intersect, we can restrict the values of $q$ to the intervals
$$D_{mn}(p,r):=\left(\dfrac{a r^{n+1} }{p^m}, \min\left(\dfrac{ r^{n+1} }{ap^m},1/36\right) \right)$$ 

We apply the Theorem \ref{main} to the family $\eS_q$ with the parameter set $D_{mn}(p,r)$ and to $S_\bj=S_3 S_1^m$ and $S_\bk=S_4 S_6^n$. We take $\br=(p,r,1/36,r,r,r)$, therefore $s_\br<1/2$ and $\bar r= 1/36$. We have $C=1$,  $C_\bk=0$ and $r_\bk=r^{n+1}$. Now since the set $K_j$ lies in the interval $[a,1]$, for $x\in K_j$ and $q',q\in D_{mn}(p,r)$ we have $|S_{\bj,q'}(x)-S_{\bj,q}(x)|=|q'-q|p^mx\ge |q'-q|p^ma$, so $c_\bj=p^m/3$. Notice also that $r^{n+1}<3p^mq$. Therefore 
$$c_\bj-C_\bk-\dfrac{r_\bj+r_\bk}{1-\bar r}>p^m\left(\dfrac{1}{3}-\dfrac{1}{35}-\dfrac{3}{35}\right)> \dfrac{p^m}{4}$$
Therefore  the set $\Da_{mn}(p,r)=\{q:S_3 S_1^m (K\mmm K_1)\cap S_4 S_6^n (K\mmm K_6)$ has the dimension less than $2s_\br$.
The same is true for the set $\Da(p,r)$ which is a countable union of the sets $\Da_{mn}(p,r)$.

This shows that\\ {\em if  $p,r\in(0,1/36)$  and $\dfrac{\log p}{\log r}$ is irrational then for Lebesgue almost all $q\in (0,1/36)$ the system $\eS$ has totally disconnected attractor with unique one-point intersection and at the same time it does not satisfy weak separation property.}

The reader may see that the properties similar to The properties {\bf 1. 2. 3.} in the previous subsection are also valid for
the systems, described above.


\begin{thebibliography}{99}

\bibitem{SSS7} {Ch.~Bandt, S.~Graf}, {\it Self-similar sets 7. A characterization of self-similar fractals with positive Hausdorff measure,}  Proc. Amer. Math. Soc. {\bf 114}:4 (1992),  995--1001.  MR1100644


\bibitem{Bar2} {B.~Barany}, Iterated function systems with non-distinct fixed points // J. Appl. Math. Anal. Appl. 383:1 (2011), pp.~244--258.

\bibitem{Bur} {N.~Bourbaki}, General Topology: Chapters 1–4 // Springer-Verlag Berlin Heidelberg, 1987.

\bibitem{Edgar} {G.~Edgar}, Measure, Topology, and Fractal Geometry //
2 ed, Springer-Verlag, New York, 2008, 272 p.

\bibitem{Fal} {K.~J.~Falconer}, Fractal geometry: mathematical foundations and applications //  J.~Wiley and Sons, New York, 1990.

 \bibitem{Hut} {J.~Hutchinson}, Fractals and self-similarity // Indiana Univ. Math. J. 30:5 (1981), pp.~713--747.

\bibitem{KT} {K.~G.~Kamalutdinov, A.~V.~Tetenov}, {Twofold Cantor sets in $\rr$}, Siberian Electronic Mathematical Reports 15 (2018), pp.~801--814.


\bibitem{Kur} {K.~Kuratowski}, {\it Topology}, vol. 1, PWN and Acad. Press, (1966). MR0217751


\bibitem{Lau} {K.~S. Lau and S.~M. Ngai}, {\it Multifractal measures and a weak separation condition, } Adv. Math. {\bf 141} (1999),  45--96. MR1667146

\bibitem{SimSolUrb} {K.~Simon, B.~Solomyak, M.~Urba\'nski}, Hausdorff dimension of limit sets for parabolic IFS with overlaps, Pacific J. Math. 201:2 (2001), pp.~441--478.


\bibitem{TKV} {A.~Tetenov, K.~Kamalutdinov, D.~Vaulin}, {Self-similar Jordan arcs which do not satisfy OSC}, arXiv:1512.00290 (2015).

\bibitem{Zer}   {M.~P.~W.~Zerner}, Weak separation properties for self-similar sets // Proc. Amer. Math. Soc. 124:11 (1996), pp.~3529--3539.



\end{thebibliography}
\end{document}